\documentstyle[12pt]{amsart}

\newcommand{\AMSLaTeX}{\protect\AmS-\protect\LaTeX}

\numberwithin{equation}{section}

\renewcommand{\thesection}{\arabic{section}}
\renewcommand{\theequation}{\thesection.\arabic{equation}}

\newtheorem{theorem}[equation]{Theorem}
\newtheorem{lemma}[equation]{Lemma}

\newtheorem{proposition}[equation]{Proposition}

\theoremstyle{definition}
\newtheorem{definition}[equation]{Definition}
\newtheorem{example}[equation]{Example}

\theoremstyle{remark}

\newcounter{itemno}

\newcommand{\twoheadrightarrow}{\rightarrow}
\def\semidirect{\Bbb o}
\newcommand{\MW}{\operatorname{MW}}
\newcommand{\SL}{\operatorname{SL}}
\newcommand{\Div}{\operatorname{Div}}

\newcommand{\dv}{\operatorname{div}}

\begin{document}
\title{Torsion Sections of Elliptic Surfaces}
\author{Rick Miranda}
\address{Dept. of Mathematics\\Colorado State University\\Ft.
Collins, CO 80523}
\email{miranda@@riemann.math.colostate.edu}
\thanks{Research supported in part by the NSF under grant
DMS-9104058 and the NSA under grant MDA 904-92-H-3022}
\thanks{\today. Set in type by \AMSLaTeX Version 1.0}
\author{Peter  Stiller}
\address{Dept. of Mathematics\\Texas A\&M University\\College
Station, TX 77843-3368}
\email{stiller@@alggeo.tamu.edu}

\maketitle

\section{Introduction}

In this article we discuss torsion sections of semistable elliptic surfaces
defined over the complex field $\Bbb C$.
Recall that a semistable elliptic surface is a fibration
$\pi:X \twoheadrightarrow C$,
where $X$ is a smooth compact surface,
$C$ is a smooth curve,
the general fiber of $\pi$ is a smooth curve of genus one,
and all the singular fibers of $\pi$ are semistable,
that is, all are  of type $I_m$ in Kodaira's notation (see
[K]).
In addition, we assume that the fibration $\pi$ enjoys a section
$S_0$;
this section defines a	zero for a group law in each fiber, making the
general fiber an elliptic curve over $\Bbb C$.

The Mordell-Weil group of $X$, denoted by $\MW(X)$,
is the set of all sections of $\pi$,
which is known to form a group under fiber-wise addition;
the section $S_0$ is the identity of $\MW(X)$.
Note that any section $S \in \MW(X)$ meets one and only one component
of each fiber.

Now given any section $S \in \MW(X)$,
one can ask the following two questions.
First, which components of the singular fibers of $X$ does $S$ meet,
and
second, exactly where in these components does $S$ meet them?

When $S$ is a torsion section,
the first question was addressed by Miranda in
	\cite{miranda2}.
To describe those results we require some notation.
Recall that a singular semistable fiber of
 type $I_m$ is a cycle of $m$ ${\Bbb P}^1$'s.
Suppose that our elliptic fibration $\pi:X \to C$
has $s$ such singular fibers $F_1,\dots,F_s$,
with $F_j$ of type $I_{m_j}$.
Choose an ``orientation'' of each fiber $F_j$ and write the $m_j$
components
of $F_j$ as
\[
C^{(j)}_0,C^{(j)}_1,\dots,C^{(j)}_{m_j-1}.
\]
where the zero section $S_0$ meets only $C^{(j)}_0$
and where for each $k$, $C^{(j)}_k$ meets only $C^{(j)}_{k\pm 1} \mod
m_j$.
If $m_j = 1$, then $F_j = C^{(j)}_0$ is a nodal rational curve
with self-intersection $0$.
If $m_j \geq 2$, then each $C^{(j)}_k$ is a smooth rational curve
with self-intersection $-2$.

Given a section $S$ of $X$,
and orientations of all of the singular fibers $F_j$,
inducing a labeling of the components as above,
we define ``component numbers'' $k_j(S)$
to be the index of the component in the $j^{th}$ fiber $F_j$
which a given section $S$ meets.  That is,
\[
S \text{ meets } C^{(j)}_{k_j(S)} \text{ in the fiber } F_j.
\]
This assigns to each section $S$, an $s$-tuple
$(k_1(S),\dots,k_s(S))$.
The component number $k_j$ can be taken to be defined $\mod m_j$
once the orientation of $F_j$ is chosen;
if the orientation is unknown,
then $k_j$ is defined  $\mod m_j$ only up to sign.
Note that $k_j(S_0) = 0$ for every $j$;
indeed, after choosing orientations,
the assignment of component numbers can be considered a group
homomorphism
from $\MW(X)$ to $\bigoplus\limits_j \Bbb Z / m_j\Bbb Z$.

Now suppose that $S$ is a torsion section of prime order $p$.
In this case, since $p\cdot S = S_0$,
we must have $p k_j(S) \equiv 0 \mod m_j$ for every $j$.
If $p$ does not divide $m_j$, then this forces $k_j(S) = 0$.
However if $m_j = p n_j$,
then $k_j(S)$ can a priori be any one of the numbers
$i n_j$, for $i = 0,\dots, p-1$.
This multiple $i$ measures, in some sense, how far around the cycle
$S$ is from the zero-section $S_0$ in the $j^{th}$ fiber.

Of course, changing the orientation in the fiber $F_j$
will have the effect of changing this multiple $i$ to $p-i$
(if $i \neq 0$; $i=0$ remains  unchanged).
We are therefore led to a definition of the following quantities.
Let $M_i(S)$ denote the number of singular fibers
where $k_j(S) = i n_j$ or $k_j(S) = (p-i) n_j$
(weighted by the number $m_j$ of components in the fiber $F_j$),
and then divided by the total number $\sum_j m_j$ of components:
\[
M_i(S) = { \sum\limits_{j \text{ with } k_j(S) = i n_j \text{ or }
(p-i)n_j} m_j \over \sum\limits_j m_j }.
\]
We may view  $M_i$  as a probability,
since it is the fraction of the fibers where $S$ meets ``distance
$i$''
away from the zero section.

The main result of \cite{miranda2} is the following:

\begin{theorem}
If $S$ is a torsion section of odd prime order $p$,
then $M_i(S) = 2p/(p^2-1)$ if $i \neq 0$,
and $M_0(S) = 1/(p+1)$.
\end{theorem}

Notice that, firstly, these fractions are independent of $S$,
and secondly, that they are independent of $i$ (for $i$ non-zero).
Thus we obtain an ``equidistribution'' property for these component
numbers.
The proof of the above theorem given in \cite{miranda2}
used only basic facts about elliptic surfaces and some intersection
theory.

In this paper we will address
the second of the two questions raised above
in a similar spirit.
Again let $S$ be a torsion section of order $p$,
and consider those fibers where $k_j(S) = 0$,
i.e., where $S$ and $S_0$ meet the same component.
By the above result, this happens exactly $1/(p+1)$ of the time,
where each fiber is weighted by the number $m_j$ of components it has.
Each identity component $C_0^{(j)}$
may be naturally identified as a group with $\Bbb C^*$,
by sending the two points of intersection with the neighboring components
to $0$ and $\infty$, and the point where $S_0$ meets $C_0^{(j)}$ to $1$.
(This identification can be made in exactly two ways,
corresponding to which node is sent to $0$ and which to $\infty$.)
Given such a coordinate choice on $C_0^{(j)}$,
the section $S$ of order $p$ will hit $C_0^{(j)}$ at a point
whose coordinate is a $p^{th}$ root of unity.

Let us denote by $\ell_j(S)$ that integer in $[1,p-1]$
such that $S \cap C_0^{(j)}$ has coordinate\break $\exp(2\pi i
\ell_j(S)/p)$.
(Since torsion sections can never meet, $\ell_j(S)$ cannot equal
$0$.)
Note that $\ell_j(S)$ is defined only when $k_j(S)=0$
and may be thought of as being defined modulo $p$;
in addition, if we switch the roles of $0$ and $\infty$,
we see that $\ell_j(S)$ is replaced by $p - \ell_j(S)$.
Thus, combinatorially, we are in a situation identical
to the one for the component numbers $k_j$.
We will call these numbers $\ell_j(S)$ {\em root of unity numbers}
for the section $S$.

Now define numbers $R_i(S)$ to be the total number of fibers
having $k_j(S)=0$ and having $\ell_j(S) = \pm i \mod p$
(weighted by the number of components $m_j$), and then
divided by the total (weighted) number of fibers with $k_j(S)=0$:
\[
R_i(S) = { \sum\limits_{j \text{ with } k_j(S)=0 \text{ and }
\ell_j(S) =
\pm i} m_j \over \sum\limits_{j \text{ with } k_j(S)=0} m_j }.
\]
The main theorem of this article is an equidistribution property
for these fractions $R_i(S)$:
namely, if $p$ is odd, $R_i(S) = 2/(p-1)$,
independent of $S$ and $i$.

These results rely on computations on elliptic modular surfaces.
We also give as a by-product
the equidistribution for the component numbers $k_j(S)$,
both by an explicit computation
and via a relationship between
the component numbers and the root of unity numbers.
Finally we develop an Abel theorem
for a singular semistable elliptic curve
and use it to compute a Weil pairing
on the singular semistable fibers of an elliptic surface.
Using this Weil pairing, we may also discover the duality between the
component numbers and root of unity numbers.

\section{Equidistribution for the roots of unity via the universal
property}

In this section we will give a proof of the following theorem.

\begin{theorem}
\label{main_theorem}
Fix a prime number $p$, and let $S$ be a section
of a smooth semistable elliptic surface which is torsion of order
$p$.
Then
\[
R_1(S) = 1 \text{ if }p=2, \text{ and }
R_i(S) = 2/(p-1) \text{ if } p \text{ is odd.}
\]
\end{theorem}

\begin{pf}
First note that if $p=2$, then the only possible value for $\ell_j$
is one,
so that certainly $R_1(S) = 1$.  Similarly, if $p=3$, the only values
for
$\ell_j$ are $1$ or $2$, which are inverse mod $3$; therefore $R_1(S)
= 1$
again.	Thus we may assume that $p$ is an odd prime at least five.

Second, note that it is immediate to check that if $S$ is a section
of
$\pi:X \to
C$, and if $F:D \to C$ is an onto map of smooth curves, then $S$
induces
a section $S'$ on the pull-back surface $\pi':X' \to D$, which is
also
torsion of order $p$. Moreover, it is easy to see that $R_i(S') =
R_i(S)$
for every $i$.

If $p \geq 5$, then we may consider the elliptic modular surface
$Y_1(p)$ over the modular curve $X_1(p)$,
which is defined using the congruence subgroup $\Gamma_1(p)$ of
$SL(2,\Bbb
Z)$ given by

\[
\Gamma_1(p) = \{
\begin{pmatrix} a & b \\ c & d \end{pmatrix}
\;|\; a,d \equiv 1 \mod p, c \equiv 0 \mod p \}.
\]
Note that $X_1(p)$ is a fine moduli space for elliptic curves
with a torsion section of order $p$,
and that $Y_1(p) \to X_1(p)$ is the universal family
(see \cite{shimura}, [Shd1] or \cite{cox-parry}).
This elliptic modular surface also has a natural section $T$ of order $p$,
and every elliptic surface with a section $S$ of order $p$
may be obtained via pull-back from this modular surface in
such a way that $S$ is the pullback of $T$.
By the above remark concerning the constancy
of the fraction $R_i$ under pull-back,
it suffices to prove that $R_i(T) = 2/(p-1)$,
for each $i = 1,\dots, (p-1)/2$;
in other words, we need only verify the statement of the theorem
for the universal section $T$ of the modular surface $Y_1(p)$ over
the
modular curve $X_1(p)$.

Now the elliptic modular surface has exactly $p-1$ singular fibers,
occuring over the cusps of the modular curve $X_1(p)$,
of which half are of type $I_1$ and half of type $I_p$ depending on
the
order of the cusp:\ the total ``weight'' of the singular fibers,
which is equal to the Euler number of the modular surface,
is $\sum_j m_j = (p^2-1)/2$.
For the $I_1$ fibers, the component numbers $k_j(T)$ must be $0$,
contributing $[(p-1)/2] / [(p^2-1)/2] = 1/(p+1)$ to the weighted
fraction
$M_0(T)$.
However, by the results of \cite{miranda2}  mentioned in the
Introduction,
the total weighted fraction $M_0(T)$ is $1/(p+1)$,
and therefore the component number must be non-zero for all the
fibers
of type $I_p$.

At this point, label the $I_1$ fibers as $F_1,\dots, F_{(p-1)/2}$.
Choose an ``orientation'' of each $I_1$ fiber,
which in essence means give an isomorphism of its smooth part
with $\Bbb C^*$, such that the zero section $T_0$ meets $F_j$
in the point corresponding to the number $1$.

Denote by $T_1=T,T_2=2T,\dots,T_{p-1}=(p-1)T$
the non-zero multiples of the universal section $T$.
Note that $T_i$ and $T_{p-i}$ meet each $I_1$ fiber $F_j$
in inverse roots of unity;
in other words, using the notation of the Introduction,
$\ell_j(T_i) = - \ell_j(T_{p-i})$ mod $p$ for every $i$ and $j$.
Construct a square matrix $Z$ of size $(p-1)/2$
whose ${ij}^{th}$ entry is the pair of integers
$\pm \ell_j(T_i) = \{\ell_j(T_i),\ell_j(T_{p-i})\}$.

Now two torsion sections of an elliptic surface never meet,
so the sections $\{T_i\}$ are all disjoint.
In particular, if we fix an index $j$,
in the $j^{th}$ column of this matrix $Z$,
we must have $p-1$ distinct integers;
since the integers come from $[1,p-1]$,
each integer in this range occurs exactly once in each column of $Z$.

Next fix an integer $i \in [1,(p-1)/2]$.
By the universality of the modular surface,
there is an automorphism of the surface fixing the zero section
and carrying $T_i$ to $T = T_1$.
Therefore, the entries along the $i^{th}$ row must be the same as the
entries
along the first row, but permuted in some way
(indeed, permuted as the automorphism permutes the $I_1$ fibers
$F_j$).

As noticed above, each integer in $[1,p-1]$
appears exactly once in the first column of $Z$.
Hence each integer appears in some row of $Z$.
By the above remark, each integer must then appear in the first row,
and indeed in every row.
This is exactly the statement that $R_i(T) = 2/(p-1)$ in this case.
\end{pf}

\section{Equidistribution for the roots of unity via explicit
computation}

In this section we will re-prove Theorem \ref{main_theorem}
via an explicit computation.
Fix an odd prime $p \geq 5$,
and let $\Gamma = \Gamma_1(p)$ be the relevant modular group
defined in the proof of the Theorem.
Let $\Bbb H$ denote the upper half-plane.
Form the semi-direct product
$\tilde{\Gamma} = \Gamma \semidirect \Bbb Z^2$
and recall that $\tilde{\Gamma}$ acts on $\Bbb H \times \Bbb C$ by
\[
(\begin{pmatrix} a & b \\ c & d \end{pmatrix},(m,n))\cdot (\tau,w) =
({a\tau + b \over c\tau + d},{1 \over c\tau + d}(w + m\tau + n)).
\]
Denote the quotient $\Bbb H \times \Bbb C / \tilde{\Gamma}$ by
$Y^0_1(p)$,
and the quotient $\Bbb H / \Gamma$ by $X^0_1(p)$.
The natural map $\pi:Y^0_1(p) \to X^0_1(p)$
induced by sending $(\tau,w)$ to $\tau$
is a smooth elliptic fibration.

The curve $X^0_1(p)$ is not compact,
but these quotients and the fibration $\pi$
extend to the natural compactifications
$\pi:Y_1(p) \to X_1(p)$.
The compact curve $X_1(p)$ is obtained by adding $p-1$ points
(called cusps) to the open curve $X^0_1(p)$,
over which the elliptic fibration has singular curves.
These cusps correspond to an equivalence class of rational points
on the boundary of the upper half-plane
and the vertical point at infinity, $i \infty$.

For this modular group $\Gamma_1(p)$,
 representatives for the cusps can be taken to be the rational
numbers
\[
{r \over p}, \text{ with }1 \leq r \leq (p-1)/2, \text{ and }
{1 \over s} \text{ with }1 \leq s \leq (p-1)/2.
\]
(See \cite{cox-parry}.)
The cusp at $i \infty$ is equivalent to the cusp $1/p$.

Over the cusps of the form $r/p$,
we have singular fibers for $\pi$ of type $I_1$;
over the cusps of the form $1/s$, we have singular fibers of type
$I_p$.

There are exactly $p$ sections of the map $\pi$,
which are induced by letting
\[
w(\tau) = {\alpha \over p}
\]
for $0 \leq \alpha \leq p-1$.
(See \cite{shioda}.)
The zero-section $T_0$ for $\pi$ is of course defined by $w(\tau) = 0$,
while the universal section $T$ can be defined by $w(\tau) = 1/p$.

A local coordinate for the modular curve $X_1(p)$ about the point
corresponding to the cusp at $\tau = i \infty$
is $q_\infty = \exp(2\pi i \tau)$;
the fibers of the modular surface itself
may be locally represented near this point
as $\Bbb C / (\Bbb Z + \Bbb Z\tau) \cong \Bbb C^* / q_\infty^{\Bbb Z}$.
Therefore, if we choose $b,d \in \Bbb Z$ such that $rd-bp=1$,
a local coordinate about the point corresponding to the cusp at $\tau = r/p$
is $q_r = \exp(2\pi i (d\tau - b)/(-p\tau + r) )$.
This essentially is a translation
(via the element
$\gamma = \begin{pmatrix} d & -b \\ -p & r \end{pmatrix}$ of
$\Gamma_0(p) \subset \SL(2, \Bbb Z)$)
from the coordinate system at $\tau = r/p$ to that at $\tau = i\infty$.

Notice that when we make this change of coordinates,
we can consider the fibers near the cusp $\tau = r/p$
again as $\Bbb C^* / q_r^{\Bbb Z}$,
where this time we have $q_r = \exp(2\pi i (d\tau - b)/(-p\tau + r) )$;
in other words, locally we are taking the complex plane and dividing
by the lattice generated by $1$ and $\gamma\tau$.
This is, up to a homothety factor $-p\tau + r$,
the lattice generated by $-p\tau + r$ and $d\tau - b$,
which is the original lattice in the complex plane
with the original variable $w$.

Now consider the section defined by $w(\tau) = \alpha/p$;
in terms of this new basis for the lattice, we have
$w(\tau) = \alpha(d\tau - b) + (\alpha d/p)(-p\tau + r)$.
Thus after applying the homothety,
and obtaining a new variable $w' = w/(-p\tau + r)$ in the complex
plane,
we see that this section is definded by
$w'(\tau) = (\alpha d/p) + \alpha(\gamma \tau)$.

After exponentiating and letting $\tau$ approach $r/p$,
we see that the section approaches the point $\exp(2\pi i\alpha
d/p)$.
This is the root of unity which we have desired to compute. Note that
had
we used $-\gamma$ instead of $\gamma$ we would have gotten the
inverse root
of unity and that $d$ is determined only $\bmod\ p$.

Using the notation of the Introduction,
let us label the $I_1$ fibers of the modular surface
as $F_1,\dots,F_{(p-1)/2}$,
with $F_r$ corresponding to the cusp at $\tau = r/p$.
We have shown that for the section $T_\alpha = \alpha T$,
\[
\ell_r(T_\alpha) = \alpha r^{-1},
\]
where the value of $r^{-1}$ is taken modulo $p$.

This gives an alternate proof of Theorem \ref{main_theorem},
since for any fixed $\alpha$, these values are equidistributed
among the nonzero classes mod $p$, up to sign.

\section{Equidistribution for the component numbers via explicit
computation}

The main theorem of \cite{miranda2}
concerning the equidistribution of the component numbers
can also be proved by appealing
to the universal property of the modular surface
$\pi:Y_1(p) \to X_1(p)$.
Using the notation of the Introduction, one sees that
the fractions $M_i(S)$ for a torsion section $S$ of order $p$
remain unchanged under base change. Thus
 it suffices to check that they have the correct values,
namely those given by Theorem~1.1,
\ for the universal section $T$ of the modular surface.
This we do in this section.

As noted in the previous section,
the modular surface contains the fibers
$F_1,\dots,F_r,\dots,\break F_{(p-1)/2}$
of type $I_1$ over the cusps represented by the rational numbers
$r/p$, for $1 \leq r \leq (p-1)/2$.
For these fibers, since they have only one component,
the component number $k_r(T)$ is zero,
contributing $[(p-1)/2]/[(p^2-1)/2] = 1/(p+1)$ to the fraction
$M_0(T)$.

The rest of the singular fibers of the modular surface
will be denoted by $F_{(p+1)/2},\dots,F_{p-1}$ where
 $F_{p-s}$ will be taken as the fiber over the cusp represented by
the rational number $1/s$, for $1 \leq s \leq (p-1)/2$.

The component number theorem then follows from the next computation.

\begin{proposition}
\label{component_number_prop}
With the above notation, if $1 \leq s \leq (p-1)/2$, then
$k_{p-s}(T) = \pm s$
(the indeterminacy of the sign being due to
the choice of orientation of the singular fiber $F_{p-s}$).
\end{proposition}

\begin{pf}
It is convenient to again transport the computation to $\tau =
i\infty$
as was done in the previous section.
Fix an $s$, with $1 \leq s \leq (p-1)/2$,
and let $\gamma = \begin{pmatrix} 1 & 0 \\ -s & 1 \end{pmatrix}$.
As $\tau$ approaches $1/s$, $\gamma\tau = \tau/(1-s\tau)$ approaches
$i \infty$.
Since the fiber $F_{p-s}$ is of type $I_p$,
the local model for the fibers of the modular surface near this cusp
can be taken to be
$\Bbb C^* / q^{p\Bbb Z}$,
where $q = \exp(2\pi i [\tau/p(1-s\tau)]) = \exp(2\pi i
\gamma\tau/p)$
is the local coordinate on the modular curve $X_1(p)$ near this cusp.

The new generators for the lattice with these coordinates
are $\tau, (1-s\tau)$ instead of $\tau, 1$.
The universal section $T$, which is defined by $w(\tau) = 1/p$,
is given in terms of these generators by the formula
$w(\tau) = (1/p)(1-s\tau) + (s/p)\tau$.
In terms of the new variable $w' = w/(1-s\tau)$,
we have $T$ defined by the formula
$w'(\tau) = 1/p + s(\gamma\tau/p)$.

Exponentiating this so as to use the local coordinate $q$,
we see that the universal section $T$ is defined in the general fiber
$\Bbb C^* / q^{p\Bbb Z}$ by
$\exp(2\pi i [1/p + s(\gamma\tau)/p]) = \zeta_p q^s$,
where $\zeta_p = e^{2\pi i/p}$.

Now using the toric description of the surface near this cusp,
(see \cite[Chapter One, Section 4]{AMRT}),
the exponent of $q$ (modulo $p$) in this formula
governs which component is hit by the section.
In particular, the universal section hits component $s$
in the fiber over the cusp represented by $1/s$. Note that using
$-\gamma$
instead of $\gamma$ would result in $-s$ instead of $s$. \end{pf}

\section{The Canonical Involution}

The matrix $A = \begin{pmatrix} 0 & -1 \\ p & 0 \end{pmatrix}$
normalizes the congruence subgroup $\Gamma_1(p)$,
and therefore induces an automorphism of the modular curve $X_1(p)$.
Since $A^2$ is a constant matrix, this automorphism
(which we will also call $A$)
is an involution, called the {\em canonical} involution, which
in terms of the variable $\tau$,  takes $\tau$ to $-1/p\tau$.
The involution $A$ also  permutes the cusps,
exactly switching the two sets of $(p-1)/2$ cusps;
indeed, the rational number $r/p$ is taken by $A$ to the number $-1/r$,
which is equivalent to the number $1/r$ under the action of $\Gamma_1(p)$.
Thus the $I_1$ cusp represented by $r/p$
is switched via $A$ with the $I_p$ cusp represented by $1/r$.

We want to simply remark that under this correspondence,
the $I_1$ cusp having $u$ as the root-of-unity number
for the universal section $T$
is paired with the $I_p$ cusp having $u^{-1}$
as the component number for $T$.

Note that for the $I_1$ fibers,
all sections have component number $k_j$ equal to zero,
and all non-zero sections have root-of-unity number $\ell_j$ in
$G(p) = {(\Bbb Z/p)}^\times/\pm 1$.
Moreover, for the $I_p$ fibers,
all non-zero sections have component number $k_j$
also in this value group $G(p)$.
For a cusp $x$ of type $I_1$,
denote the root-of-unity number in $G(p)$ by $\ell_x$;
for a cusp $x$ of type $I_p$,
denote the component number in $G(p)$ by $k_x$.
Then the above statement can be re-phrased as follows.
\begin{equation}
\label{ALcuspformula}
\text{For every cusp } x \text{ of type } I_1,\;\;\;
\ell_x = k_{Ax}^{-1} \text{ in } G(p).
\end{equation}

This is a manifestation of a certain duality
between the two notions.
This we will explore further in the next sections,
via a version of the Weil pairing on singular fibers
of elliptic surfaces.
Before proceeding to this, we want to make some elementary remarks
concerning the modular surface and this involution.

As above, denote the modular surface over $X_1(p)$ by $\pi:Y_1(p) \to
X_1(p)$.
Let $\pi':Y_1(p)' \to X_1(p)$ denote
the pull-back of $\pi$ via the involution $A$.
This operation exactly switches the fibers over the cusps,
so that the fiber of $\pi'$ over an $r/p$ cusp is now of type $I_p$,
and the fiber of $\pi'$ over a $1/s$ cusp is of type $I_1$.
The universal section $T$ of $\pi$ pulls back to a section $T'$ of $\pi'$.
(We note that this surface is the modular surface for the group $\Gamma^1(p)$.)

An alternate way of constructing this elliptic fibration $\pi'$
is to take the original modular surface $\pi:Y_1(p) \to X_1(p)$
and divide it, fiber-by-fiber, by the subgroup generated by the
universal torsion section $T$.
Over the cusps represented by $r/p$,
the original modular surface has $I_1$ fibers,
which are rational nodal curves;
in the quotient there is again a rational nodal curve,
but the surface acquires a rational double point of type $A_{p-1}$
at the node, and the minimal resolution of singularities produces
a fiber of type $I_p$.
Over the cusps represented by $1/p$,
we have $I_p$ fibers in the original surface;
in the quotient the $p$ components are all identified
to a single $I_1$ component.
For details concerning this quotient construction,
the reader may consult \cite{miranda-persson3}.

The section $T'$ in this view comes not from the original section
$T$,
but from a $p$-section of the modular surface,
consisting of a coset of the cyclic subgroup generated by $T$
in the general fiber of $\pi$.
(The universal section $T$ and all of its multiples
of course descends to the zero-section of $\pi'$.)
The formula (\ref{ALcuspformula}) implies that
\[
k_x(T') = \ell_x(T)^{-1}
\]
in the group $G(p)$.

\section{The function group of a semistable elliptic curve}

In the next section we will develop a version of the Weil pairing
on a semistable singular fiber of an elliptic surface
(that is, a fiber of type $I_m$).
The essential ingredient in the definition of the Weil pairing
on a smooth elliptic curve is Abel's theorem;
see for example \cite[Chapter III, Section 8]{silverman}.
Therefore we require a version of Abel's theorem
for a fiber of type $I_m$,
and this in turn requires a notion of appropriate rational functions
on the degnerate fiber.
In this section we describe the set of functions which we will use;
Abel's theorem is an immediate consequence of the definition.

For a nodal fiber $F$ of type $I_1$,
which is irreducible,
we have the function field of rational functions on $F$,
which may be identified with the field of rational functions on $\Bbb
P^1$
which take on the same value at $0$ and $\infty$.
The non-zero elements of this field form a multiplicative group,
which has as a subgroup ${\cal K}$ the group of functions which are
regular and nonzero at the node.
This group $\cal K$ comes equipped with a divisor function to the
group
of divisors $\Div(F^{sm})$ supported on the smooth part $F^{sm}$ of
$F$.
Since $F^{sm}$ is a group,
there is a natural map $\Phi:\Div(F^{sm}) \to F^{sm}$
which takes a formal sum of points of $F^{sm}$ to the actual sum in
the group.
It is easy to check that Abel's theorem holds:
a divisor $D \in \Div(F^{sm})$ is the divisor of a function $f \in
{\cal K}$
if and only if $\deg(D) = 0$ and $\Phi(D) = 0$ in the group law of
$F^{sm}$.

We will now extend this to fibers of type $I_m$ with $m \geq 2$.
With this assumption there is no field of functions to employ,
since the fiber $F$ is no longer irreducible.
Thus we must find a multiplicative group of appropriate functions
without the aid of an associated field of rational functions.

For this we rely on the existence of a set of coordinates on the
components
of $F$, which are adapted nicely to both the group law on $F^{sm}$
and to the elliptic surface on which $F$ lies.

\begin{definition}
Suppose that the singular fiber $F$ (which is assumed to be of type
$I_m$)
has components $C_0,C_1,\dots,C_{m-1}$,
with the zero section $S_0$ meeting $C_0$
and $C_j$ meeting only $C_{j\pm 1}$ for each $j$.
Let $u_j$ be an affine coordinate on $C_j$ for each $j$.
The set $\{u_j\}$ of coordinates will be called
a {\em standard set of affine coordinates on $F$}
if the following conditions hold:
\begin{itemize}
\item[a)] For each $j$,
$u_j = 0$ at $C_j \cap C_{j-1}$ and $u_j = \infty$ at $C_j \cap
C_{j+1}$.
\item[b)] The map $\alpha:\Bbb C^* \times \Bbb Z/m \to F^{sm}$
sending a pair $(t,j)$ to the point $u_j=t$ in component $C_j$ of $F$
is an isomorphism of groups.
\item[c)] For each $j$, if we set $v_j = 1/u_j$
to be the affine coordinate on $C_j$ near $u_j = \infty$,
then $v_j$ and $u_{j+1}$ extend to coordinate functions $V_j$ and
$U_{j+1}$
on the elliptic surface near the point $C_j \cap C_{j+1}$,
such that $V_j = 0$ defines $C_{j+1}$ and $U_{j+1} = 0$ defines $C_j$
near this point.
\end{itemize}
\end{definition}

\begin{proposition}
Let $F$ be a fiber of type $I_m$ on a smooth elliptic surface $X$.
Then a standard set of affine coordinates exist on $F$.
Moreover, given the ordering of the components,
there are exactly $m$ such standard sets of affine coordinates on
$F$.
\end{proposition}

\begin{pf}
The existence of a standard set of affine coordinates on $F$
is an immediate consequence of the local toric description
of a smooth semistable elliptic surface near a singular fiber of type
$I_m$
given in \cite[Chapter One, Section 4]{AMRT}.  Indeed, the standard
set
of affine coordinates is exactly the set of coordinates on the toric
cover
described there, which descend nicely to $X$.

Note that by $a)$ and $b)$, the coordinate $u_0$ on $C_0$ is
determined:
it must be $0$ at $C_0\cap C_{m-1}$, $\infty$ at $C_0\cap C_1$,
and $1$ at $S_0\cap C_0$.
Similarly, each coordinate $u_j$ is determined by the point $u_j=1$,
which must be one of the  points of $C_j$ of appropriate order.
Thus there are exactly $m$ possibilities for $u_1$.
By $b)$, once $u_1$ is determined, so is $u_j$.
It is easy to check that these $m$ different possibilities
all give standard sets (once it is known that one of them does).
\end{pf}

Note that if $\{u_0,\dots,u_{m-1}\}$ is a standard set of affine
coordinates
on $F$,
then any other standard set is of the form
$\{u_j' = \zeta^j u_j\}$
for some $m^{th}$ root of unity $\zeta$.

We can now define the group of functions $\cal K$
which plays the role of rational functions on $F$.
A bit of notation is useful.
Suppose that $g(u)$ is a nonzero  rational function of $u$.
Firstly, define $n_0(g)$ to be the order of $g$ at $u=0$,
and $n_\infty(g)$ to be the order of $g$ at $u=\infty$
(which is the order of $g(1/v)$ at $v=0$).
These integers are of course independent
of the choice of affine coordinate $u$.
Moreover, we have that $g(u)/u^{n_0(g)}$ has a finite value at $u=0$,
and $g(u)u^{n_\infty(g)}$ has a finite value at $u=\infty$.
Secondly, define $c_0(g)$ to be
the value of $g(u)/u^{n_0(g)}$ at $u=0$,
and define $c_\infty(g)$ to be
the value of $g(u)u^{n_\infty(g)}$ at $u=\infty$.
These constants do depend on the choice of coordinate $u$;
they are simply the lowest coefficient of the Laurent series for $g$
expanded about $u=0$ and $u=\infty$.

If we write
\begin{equation}
g(u) = \alpha u^\ell
\prod_{i=1}^e {(u - \lambda_i)} / \prod_{k=1}^f {(u - \mu_k)},
\end{equation}
with $\alpha$, $\lambda_i$ and $\mu_k$ nonzero, and $\ell \in \Bbb
Z$,
then
\begin{eqnarray}
n_0(g) &= & \ell\\
n_\infty(g) &=& f - e - \ell \\
c_0(g) &=&
\alpha {(-1)}^{e + f} \prod_i \lambda_i / \prod_k \mu_k,
\text{ and } \\
c_\infty(g) &=& \alpha.
\end{eqnarray}

Fix a standard set $\{u_j\}$ of affine coordinates on $F$.
Define $\cal K$ to be the set of $m$-tuples of nonzero rational
functions
$(g_0(u_0), g_1(u_1), \dots, g_{m-1}(u_{m-1}))$
satisfying the following conditions:
\begin{itemize}
\item[a)] For each $j$, $n_\infty(g_j) + n_0(g_{j+1}) = 0$.
\item[b)] For each $j$, $c_\infty(g_j) = c_0(g_{j+1})$.
\item[c)] $\sum_{j=0}^{m-1} n_0(g_j) = 0$.
\end{itemize}

Note that condition $a)$ says that
the functions $\{g_j\}$ have opposite orders at the nodes of $F$,
and condition $b)$ says that
with respect to the standard set of affine coordinates,
they have the same leading coefficients in their Laurent series at
these nodes.
These are local conditions about each of the nodes.
The final condition $c)$ is a global condition,
saying that the orders sum to zero upon going around the cycle of
components.
We note that $\cal K$ is a multiplicative group
(the operation being defined component-wise).

These conditions are motivated by a notion of restriction of
functions
from the surface to the fiber.
Suppose that $G$ is a rational function on $X$.
We may uniquely write its divisor, near the fiber $F$, as $\dv(G) = H
+ V$,
where $V$ is the ``vertical'' part of the divisor
consisting of linear combinations of components of $F$,
and $H$ is the ``horizontal'' part of the divisor
consisting of linear combinations of multi-sections for the fibration
map.
If one restricts to the general fiber near $F$,
one only sees the contributions from $H$.
We want to make the following regularity assumption for the function
$G$:
\begin{center}
$\{*\}$\;No curve appearing in the horizontal part $H$ passes through
any node of $F$.
\end{center}
Under this assumption, we see that the zeroes and poles of $G$,
as we approach the singular fiber $F$,
survive in the smooth part $F^{sm}$ of $F$.
Therefore we have a chance of obtaining a limiting version of Abel's
theorem.

Take then such a rational function $G$,
and let us define a ``restriction'' to the fiber $F$,
which will be an element of the group $\cal K$.
Fix a standard set of affine coordinates $\{u_j\}$ on $F$,
and also assume that we have normalized the base curve
so that $\pi = 0$ along $F$.
Write $V = \sum_j r_j C_j$ as the vertical part of $\dv(G)$.
For each $j$, consider the ratio $G/\pi^{r_j}$;
this is a rational function on $X$
which does not have a zero or pole identically along $C_j$
(since $\pi$ has a zero of order one along $F$).
Restricting this function to $C_j$ gives
a nonzero rational function $g_j(u_j)$.

\begin{lemma}
If $G$ satisfies the regularity condition $\{*\}$,
then the $m$-tuple $(g_0,\dots,g_{m-1})$ lies in $\cal K$.
\end{lemma}

\begin{pf}
Near the point $C_j \cap C_{j+1}$,
we have local coordinates $U_{j+1}$ and $V_j$ on $X$
as in the definition of a standard set of affine coordinates.
The curve $C_j$ is defined locally by $U_{j+1} = 0$,
and $C_{j+1}$ by $V_j = 0$.
Moreover the fibration map $\pi$ is locally of the form $\pi =
U_{j+1} V_j$.
Hence near $C_j \cap C_{j+1}$, we may write $G$ as
\[
G(U_{j+1},V_j) = U_{j+1}^{r_j} V_j^{r_{j+1}} L(U_{j+1},V_j)
\]
where the condition $\{*\}$ implies that $L(0,0) \neq 0$.
With this notation we have $g_j(v_j) = v_j^{r_{j+1}-r_j} L(0,v_j)$
and $g_{j+1}(u_{j+1}) = u_{j+1}^{r_j-r_{j+1}} L(u_{j+1},0)$.

Thus we see  that $n_\infty(g_j) = r_{j+1}-r_j$
and $n_0(g_{j+1}) = r_j-r_{j+1}$,
proving that condition $a)$ of the definition of $\cal K$ is
satisfied.

We also have that $c_\infty(g_j) = c_0(g_{j+1}) = L(0,0)$,
giving us condition $b)$.

Finally, the sum $\sum_j n_0(g_j)$ telescopes to $0$,
showing that condition $c)$ holds.
\end{pf}

{}From this point of view, the definition of $\cal K$,
though at first glance rather ad-hoc,
is actually quite natural.
Motivated by the above, we call $\cal K$
the {\em function group} of the singular fiber $F$.

We have a divisor map $\dv:{\cal K} \to \Div(F^{sm})$,
sending an $m$-tuple $(g_0,\dots,g_{m-1})$
to the formal sum of the zeroes and poles of each $g_j$,
throwing away any part of the divisor at the nodes.
This map is a group homomorphism.
Recall that we also have a natural summation map $\Phi$
from $\Div(F^{sm})$ to $F^{sm}$.
Abel's theorem for $F$ can be stated as follows.

\begin{theorem}
\label{abel}
A divisor $D \in \Div(F^{sm})$ is the divisor of an element of $\cal
K$
if and only if $\deg(D) = 0$ and $\Phi(D) = 0$.
\end{theorem}

\begin{pf}
Let $(\underline{g}) \in {\cal K}$, and let $D = \dv(\underline{g})$.
For each $j$, write
\[
g_j(u_j) = \alpha_j u_j^{\ell_j}
\prod_{i=1}^{e_j} {(u_j - \lambda_i^{(j)})} /
\prod_{k=1}^{f_j} {(u_j - \mu_k^{(j)})}.
\]
{}From the definition of $\cal K$, we must have
\begin{eqnarray*}
\ell_j - \ell_{j+1} &=& f_j - e_j,\\
\alpha_j &=&
\alpha_{j+1} {(-1)}^{e_{j+1} + f_{j+1}}
\prod_i \lambda_i^{(j)} / \prod_k \mu_k^{(j)}, \text{ and }\\
\sum_j \ell_j &=& 0.
\end{eqnarray*}
Summing the first set of equations over $j$, we see that
\[
\deg(D) = \sum_j (e_j - f_j) = 0.
\]
Multiplying the second set of equations over $j$, and applying the
above,
we have that
\[
\prod_j {\prod_i \lambda_i^{(j)} \over \prod_k \mu_k^{(j)} } = 1
\]
which shows that the $\Bbb C^*$ part of the group element $\Phi(D)$
is
trivial.
Finally, to show that the $\Bbb Z/m$ part of $\Phi(D)$ is trivial,
we must show that $\sum_j j (e_j - f_j) = 0 \mod m$.
Writing $e_j - f_j$ as $\ell_{j+1} - \ell_j$ using the first
equation,
we see that this sum telescopes to
\[
-\ell_1 - \ell_2 - \dots - \ell_{m-1} + (m-1)\ell_0
\]
which is $0 \mod m$ by the third equation.
This completes the proof of the necessity of the conditions on $D$.

We leave the sufficiency, which is equally elementary, to the reader.
\end{pf}

We note that the only elements $(\underline{g})$ of $\cal K$
which have $\dv(\underline{g}) = 0$
are the constant elements, where $g_j = c$ for every $j$
with $c$ being a fixed nonzero complex number.
We thus have an exact sequence
\[
0 \to \Bbb C \to {\cal K} \stackrel{\dv}{\to} \Div_0(F^{sm})
\stackrel{\Phi}{\to} F^{sm} \to 0
\]
where $\Div_0(F^{sm})$ is the group of divisors of degree $0$ on
$F^{sm}$.

\begin{example}
Suppose $F$ is a fiber of type $I_{mk}$,
with a standard set of affine coordinates $\{u_j\}$.
Fix a primitive $m^{th}$ root of unity $\zeta$, and an integer
$\alpha \in
[0,m-1]$,
and let $p$ be the point of $C_{\alpha k}$
with coordinate $u_{\alpha k} =\zeta$.
We note that $p$ is a point of order $m$ in the group law of
$F^{sm}$,
so that $D = m p - m \underline{0}$ is a divisor on $F^{sm}$
with $\deg(D) = 0$ and $\Phi(D) = \underline{0}$
(where $\underline{0}$ is the origin of the group law,
i.e., the point in $C_0$ with coordinate $u_0 = 1$).

By Abel's theorem, there is an element $(\underline{g}) \in {\cal K}$
such that $\dv(\underline{g}) = D$.
If $\alpha = 0$, this element $(\underline{g})$ is
\[
g_0 = {(u_0 - \zeta)}^m / {(u_0 - 1)}^m, \;\;\; g_j = 1
\text{ for } j \neq 0.
\]
If $1 \leq \alpha \leq m-1$, we have
\begin{eqnarray*}
g_0 &=& u_0^{\alpha}/{(u_0 - 1)}^m, \\
g_j & = & u_j^{\alpha-m} \text{ for }j = 1,\dots,\alpha k - 1, \\
g_{\alpha k} &=& {(-1)}^m u_{\alpha k}^{\alpha-m} {(u_{\alpha k} -
\zeta)}^m,
\text{ and }\\
g_j &=& {(-1)}^m u_j^{\alpha} \text{ for } j = \alpha k + 1, \dots,
mk-1.\\
\end{eqnarray*}

\end{example}

\section{The limit Weil pairing on a degenerate elliptic curve}

The Weil pairing on a smooth elliptic curve $E$
is defined as follows (see \cite[Section III.8]{silverman}).
Let $S$ and $T$ be two points of order $m$ on $E$.
Choose a rational function $g$ on $E$
with $\dv(g) = {[m]}^*(T) - {[m]}^*(\underline{0})$,
where $[m]$ denotes multiplication by $m$.
Then the Weil pairing $e_m$ on the $m$-torsion points of $E$
is defined by
\[
e_m(S,T) = g(X \oplus S) / g(X)
\]
for any $X \in E$ where both $g(X\oplus S)$ and $g(X)$ are defined
and nonzero.
(We use $\oplus$ as the group law in $E$ to avoid confusion.)
The existence of the rational function $g$
relies solely on Abel's theorem for $E$,
as does the fact that $e_m$ has values
in the group of $m^{th}$ roots of unity.

In the previous section, we developed an Abel theorem
for a singular fiber $F$ of type $I_k$ on a smooth elliptic surface,
by replacing the notion of the field of rational functions
with the limit function group $\cal K$.
This allows us to define in the same way a limit Weil pairing
on the $m$-torsion points of $F$,
which we also denote by $e_m$.
In this section we compute it.

Since the $e_m$ pairing is isotropic and skew-symmetric,
it suffices to compute $e_m(S,T)$ for generators $S$ and $T$
of the group of $m$-torsion points on $F$.

Fix an integer $m$; in order that $F$ have $m^2$ $m$-torsion points,
$F$ must be of	type $I_{mk}$ for some $k$.
Since the Weil pairing is formally invariant under base change,
we may assume $k=1$ and $F$ is of type $I_m$.
Also fix a standard set of affine coordinates $\{u_j\}$ on $F$.
Let $\zeta = \exp(2 \pi i/m)$ be a primitive $m^{th}$ root of
unity.
Let $T$ be the $m$-torsion point of $F$
which is in component $C_0$
having coordinate $u_0 = \zeta$.
Let $S$ be the $m$-torsion point of $F$
which is in component $C_1$
having coordinate $u_1 = 1$.
These points $S$ and $T$ generate the group of $m$-torsion points of
$F$.

Let $\nu = \exp(2 \pi i/m^2)$ so that $\nu^m = \zeta$.
Consider the point $T' \in F$ in component $C_0$
having coordinate $u_0 = \nu$;
note that $m T' = T$,
and ${[m]}^*(T) = \sum_R (T' \oplus R)$
where the sum is taken over all $m$-torsion points of $F$.
Similarly, ${[m]}^*(\underline{0}) = \sum_R (R)$.

The element $(\underline{g}) \in {\cal K}$
such that $\dv(\underline{g}) = {[m]}^*(T) - {[m]}^*(\underline{0})$
is defined by
\[
g_j(u_j) = \zeta^{-j} {u_j^m - \zeta \over u_j^m - 1}.
\]

Now choose an $X \in C_0$ with coordinate $u_0 = x$.
The point $X\oplus S$ is then the point in $C_1$ with coordinate $x$.
Thus
\begin{eqnarray*}
\underline{g}(X \oplus S) / \underline{g}(X) &=& g_1(x) / g_0(x) \\
&=& \zeta^{-1} {x^m - \zeta \over x^m - 1} / {x^m - \zeta \over x^m -
1} \\
&=& \zeta^{-1}.
\end{eqnarray*}
This shows that $e_m(S,T) = \zeta^{-1}$.

Let $M(t,s) = t T \oplus s S$ be the general point of order $m$ in
$F$;
note that $M(t,s)$ is in component $C_s$ and has coordinate $u_s =
\zeta^t$.
Extending the calculation above using the bilinearity and
skew-symmetry
we obtain the following.

\begin{proposition}
With the above notation, the limit Weil pairing on $F$ takes the form
\[
e_m(M(t_1,s_1),M(t_2,s_2)) = \zeta^{t_1s_2 - t_2s_1}.
\]
\end{proposition}

This limit Weil pairing,
although computed using Abel's theorem for the limit function group
$\cal K$,
is a priori dependent upon the choice
of a standard set of affine coordinates on the fiber.
However it is easy to see that	it is in fact  well-defined,
independent of this choice.
Moreover the limit Weil pairing is indeed
the limit of the usual Weil pairing
on the nearby smooth fibers of the elliptic surface.
This follows from the nature of the element $(\underline{g}) \in
{\cal K}$
used in the computation above:
each $g_j$ has degree $0$ on the component where it is defined,
and hence is the usual restriction
of a rational function $G$ on the elliptic surface.
This function $G$ can be chosen so that,
when restricted to the nearby smooth fibers,
it is the function used to define the Weil pairing there.
Therefore the limit Weil pairing is the limit of the usual Weil
pairing.

\section{The root-of-unity and component number relationship
via the Weil pairing}

Finally, we want to point out that the duality
between the root-of-unity results and the component number results,
which were mentioned in Section~4 
as being related to the canonical involution,
can be expressed also in terms of the Weil pairing.
We will compute this Weil pairing on the singular fibers
using the limit Weil pairing developed in the previous section.

First let $W$ be a torsion section of order $m$
of an elliptic surface
passing through the zero component $C_0$ of an $I_m$ fiber.
Let $\zeta = \exp(2\pi i/m)$ be a primitive $m^{th}$ root of
unity.
Assume that the $I_m$ fiber is given a standard set of affine
coordinates,
such that the point $W \cap C_0$ has coordinate $u_0 = \zeta^a$,
with $(a,m) = 1$.
In the notation of the last section, we have $W = aT$.
Let $W^*$ be the set of (local) sections $Z$ such that
$e_m(W,Z) = \zeta$.
By the computation given in Proposition~7.1,
\ we see that such a local section $Z$ is one
which passes through component $b$, where $ab \equiv 1 \mod m$.

Specializing to the case where $m$ is an odd prime $p$,
and using the root-of-unity numbers $\ell$ and the component numbers
$k$,
we see the following:
\begin{equation}
\label{W*}
\text{If } k_j(W) = 0 \text{ and } \ell_j(W) = a,
\text{ then } W^* \text{ is the set of local sections } Z \text{ with }
k_j(Z) = a^{-1}.
\end{equation}

Finally let us return to the modular surface $\pi:Y_1(p) \to X_1(p)$
and the quotient $\pi':Y_1'(p) \to X_1(p)$,
as described in Section~4.
\ In the above, set $W = T$, the universal section,
 and fix a cusp $x$ of type $I_1$ on
$X_1(p)$.
Assume that $\ell_x(T) = a$ for this cusp,
that is, $T\cap C_0$ is the point with coordinate $\zeta^a$.
The section $T'$ on the quotient is induced by the multisection of
$\pi$
which, by (\ref{ALcuspformula}),
has $k_x = a^{-1}$.
Since the Weil pairing is invariant under isogeny,
we see by (\ref{W*}) that $T'$ is exactly the image, locally,
of the set $W^*$ of sections pairing with $T$ to give value $\zeta$.
Alternatively, we may write
\[
T' = \text{ image of } {e_m(T,-)}^{-1}(\zeta).
\]

\end{document}